\documentclass[11pt]{article}
\usepackage{graphicx,mathrsfs,amssymb}
\usepackage{amsmath,latexsym, color, amsbsy,amssymb}
\usepackage{graphicx}

\makeatletter
\newcommand{\superimpose}[2]{%
	{\ooalign{$#1\@firstoftwo#2$\cr\hfil$#1\@secondoftwo#2$\hfil\cr}}}
\makeatother
\newcommand{\transversal}{\mkern-1mu\mathrel{\mathpalette\superimpose{{\top}{\scrunch{\cap}}}}\mkern-1mu}
\newcommand{\scrunch}[1]{\resizebox{\width}{.9\height}{$#1$}}

\def\ds{\displaystyle}
\def\forall{\hbox{for all}~}

\def\bfv{{\bf v}}

\def\ve{\varepsilon}
\def\n{\noindent}

\def\R{\mathbb{R}}

\def\vp{\varphi}

\def\vs{\vskip 2em}

\def\C{{\cal C}}
\def\bega{\begin{array}}
	\def\enda{\end{array}}
\def\begi{\begin{itemize}}
	\def\endi{\end{itemize}}

\def\bel{\begin{equation}\label}
	\def\eeq{\end{equation}}
\def\sqr#1#2{\vbox{\hrule height .#2pt
		\hbox{\vrule width .#2pt height #1pt \kern #1pt
			\vrule width .#2pt}\hrule height .#2pt }}

\def\square{\sqr74}
\def\endproof{\hphantom{MM}\hfill\llap{$\square$}\goodbreak}

\newtheorem{theorem}{Theorem}[section]

\newtheorem{definition}[theorem]{Definition}
\newtheorem{remark}[theorem]{Remark}
\newtheorem{lemma}[theorem]{Lemma}

\begin{document}
	\title{\bf  A sharp  quantitative estimate of  critical sets}\vs
	\author{\it Andrew Murdza and Khai T. Nguyen\\
\\
		{\small Department of Mathematics, North Carolina State University}\\
		{\small ~apmurdza@ncsu.edu, ~khai@math.ncsu.edu}
}
\maketitle
\begin{abstract}

The paper  establishes a sharp quantitative estimate   for the $(d-1)$-Hausdorff measure of the critical set  of $\mathcal{C}^1$ vector-valued functions on $\R^d$. Additionally, we prove that for a generic $\C^2$ function where ``generic" is understood in the topological sense of Baire category, the critical set has a locally finite $(d-1)$-Hausdorff measure.

%
\quad\\
		\quad\\
		{\footnotesize
		\medskip
		
		\n {\bf AMS Mathematics Subject Classification.} 46T20
		}
	\end{abstract}
	
	\section{Introduction}
	\label{sec:1}
	\setcounter{equation}{0}
The well-known Morse-Sard theorem  states that for every  $k\geq \max\{d-m+1,1\}$ and $d\geq m$, the critical set  of a $\mathcal{C}^k$ function $f:\R^d\to\R^m$, which is denoted by 
	\[
	\Sigma_f~=~\left\{x\in \R^d:\mathrm{rank}(D f(x))\leq m-1\right\}, 
	\]
 can be very large, but its image is still small in the sense of Lebesgue measure in $\R^m$. Indeed, it was shown in  \cite{M1,S} that 
	\[
	\mathcal{L}^{m}\left\{f(x): x\in \Sigma_f\right\}~=~0,
	\]
	and  the Hausdorff dimension of $f(\Sigma_f)$ can be arbitrarily close to $m$. This theorem, which holds a fundamental role in singularity theory among other fields, has several variants. A version applicable to infinite-dimensional Banach manifolds was established  in \cite{SS}. Other's related works or development can be found, for example in \cite{B,AF,P}. From a different aspect,  one observes that there are  $d$ unknown variables $(x_1,x_2,\dots, x_d)$ and single equation $\det(D f(x))=0$ to verify that $x\in \Sigma_f$.  Hence, it is quite natural  to expect a generic result on  $\mathcal{H}^{d-1}$-rectifiability of $\Sigma_f$. Using the transversality lemma \cite{RT}, we offer a simple proof for this. Roughly speaking, we  show that 
	\begin{theorem}[Generic $\mathcal{H}^{d-1}$ rectifiable] For ``nearly all" $f\in\mathcal{C}^{2}(\R^d,\R^m)$, $\Sigma_f$ has a locally finite $(d-1)$-Hausdorff measure.	
		\end{theorem}

Here ``nearly all" is meant in the topological sense of Baire category: these properties should be true on a $G_{\delta}$ set, i.e., on the intersection of countably many open dense subsets. This leads to an interesting question on the quantitative analysis of the measure of $\Sigma_f$. Namely, how small can we make this measure by an $\ve$-perturbation of $f$? To express our answer to this question more precisely, we assume that   $f\in \mathcal{C}^1([0,1]^d,\R^m)$ for simplicity and define
\[
\mathcal{N}_f(\ve)~=~\inf_{\|g-f\|_{\mathcal{C}_1}\leq \ve}\mathcal{H}^{d-1}\left(\mathcal{Z}_g\right),\qquad \ve>0,
\]
to be the smallest $(d-1)$-Hausdorff measure of the critical set of $g\in \mathcal{C}^1([0,1]^d,\R^m)$ with $\|g-f\|_{\mathcal{C}^1}\leq \ve$. In Theorem \ref{N-1}, we shall establish an upper bound  on the number $\mathcal{N}_f(\ve)$, for a general $\C^1$ function $f$ by using the concept of the modulus of continuity. This bound can be expressed in terms of the inverse of the minimal modulus of  continuity of $D f$. Specifically,  we obtain the following estimate for a class of $\mathcal{C}^{1,\alpha}$-functions for $\alpha\in (0,1]$.

\begin{theorem}\label{Ho} Let $f: [0,1]^d\to\R^m$ be a $\mathcal{C}^{1,\alpha}$-function with exponent $\alpha\in (0,1]$ and $m\leq d$. Then for every $\ve>0$ sufficiently small, it holds
\[
\mathcal{N}_f(\ve)~\leq~C_{[d,\alpha]}\cdot\Bigg(\frac{\|D f\|_{\mathcal{C}^{0,\alpha}}}{\ve}\Bigg)^{1/\alpha},
\]
where  the constant $C_{[d,\alpha]}$ depends only on $d,\alpha$ and $\|D f\|_{\mathcal{C}^{0,\alpha}}=$ is the  H\"older norm of $D f$.
\end{theorem}
The proof of Theorem \ref{Ho} relies on a direct construction using a simple lemma concerning the decomposition of a  $d$-dimensional cube, basic tools in linear algebra, and a recent result in \cite{MN0} on a sharp universal bound on the ($d-1)$-Hausdorff measure of the zeros of nontrivial multivariable polynomials. Finally, the sharpness of the upper bound in terms of the power function is proved in Theorem \ref{lwb}.

\section{Generic properties of critical sets}
\setcounter{equation}{0}
In this section, we shall provide a simple proof on a generic result for critical sets of smooth functions  by using the transversality lemma which is stated as follows:  let $g:X\to Y$ be a smooth map between  two smooth manifolds $X$ of dimension $d$ and $Y$ of dimension $m$. For any smooth submanifold $W\subset Y$, we say that the smooth function $g:X\to Y$ is transverse to $W$ and write $g\transversal W$ if for all $q\in g^{-1}(W)$ 
\bel{tran}
(dg)_q(T_qX)+T_{g(q)}(W)~=~T_{g(q)}(Y),
\eeq
with $T_zZ$ being the tangent space to a smooth manifold $Z$ at a point $z\in Z$ and $(dg)_q:T_q(X)\to T_{g(q)}(Y)$ being the derivative of $g$ at $q$.
\begin{lemma}\label{Tranv}
Given a smooth manifold $\Theta$, let $\theta\mapsto\phi^{\theta}$ be a smooth map which to each $\theta\in \Theta$ associates a smooth map $\phi^{\theta}:X\to Y$, and define $\Phi:X\times\Theta\to Y$ by setting
\[
\Phi(x,\theta)~=~\phi^{\theta}(x)\quad\forall (x,\theta)\in X\times\Theta.
\]
If $\Phi\transversal W$ then the set $\{\theta\in \Theta: \phi^{\theta}\transversal W\}$ is dense in $\Theta$.
\end{lemma}
Our  result is stated as follows.
\begin{theorem}
For every given $1\leq m\leq d<\infty$, there exists a $\mathcal{G}_{\delta}$ subset $\mathcal{M}$ of $\C^2(\R^d,\R^m)$ such that for every $f\in \mathcal{M}$, it holds
\bel{F-Hau}
\mathcal{H}^{d-1}\left(\Sigma_{f}\bigcap K\right)~<~+\infty,\qquad\forall K\subset\R^{d}~~\mathrm{compact}.
\eeq
\end{theorem}
{\bf Proof.} {\bf 1.} For every $f\in \C^2(\R^d,\R^m)$, let we define the $\mathcal{C}^1$ map $\Psi^{f}:\R^d\times S^{m-1}\to\R^m$  by setting
\bel{Phi}
\Psi^f(x,{\bf v})~=~\Pi_{m}\big(D^{T} f(x){ \bf v}\big),\qquad (x,{\bf v})\in\R^d\times S^{m-1},
\eeq
where  $S^{m-1}=\left\{x\in \R^m:|x|=1\right\}$, $D^{T} f(x)$ is the transpose of the Jacobian matrix $Df(x)\in \R^{m\times d}$, and $\Pi_{m}:\R^d\to\R^m$ be a projection such that 
\[
\Pi_{m}(w)~=~(w_{1},\cdots, w_{m}),\qquad w\in \R^d.
\]
For $\nu\geq 1$, calling $\overline{B}_{\nu}$   the closed ball centered at the origin with radius $\nu$, we consider the subset of $\C^2(\R^d,\R^m)$
\[
\mathcal{M}_{\nu}~=~\left\{f\in \mathcal{C}^2(\R^d,\R^m):\Psi^f_{|\overline{B}_{\nu}\times S^{m-1}}~~\text{is~transversal to}~\{0\}\right\}.
\]
By (\ref{Phi}) and (\ref{tran}),  $\Psi^f_{|\overline{B}_{\nu}\times S^{m-1}}$ is transversal to $0$ if 
\[
\mathrm{rank} \big(D\Psi^f(x,{\bf v})\big)~=~m,\qquad (x,{\bf v})\in\overline{B}_{\nu}\times S^{m-1},
\]
and the continuity of $D\Psi^f$  implies that $\mathcal{M}_{\nu}$ is open in $\C^2(\R^d,\R^m)$. Hence, if  $\mathcal{M}_{\nu}$ is dense in $\C^2(\R^d,\R^m)$ for all $\nu\geq 1$, then the set $\mathcal{M}\doteq \ds\bigcap_{\nu\geq 1}\mathcal{M}_{\nu}$ is a $\mathcal{G}_{\delta}$ subset of $\C^2(\R^d,\R^m)$ such that  for every $f\in \mathcal{M}$, the map $\R^d\times S^{m-1} \ni (x,{\bf v})\mapsto \Psi^f(x,{\bf v})$ is  transverse to $\{0\}$. By  the implicit function theorem, the set 
\[
\mathcal{Z}_{\Psi^f}~=~\left\{(x,{\bf v})\in \R^d\times S^{m-1}:\Psi^f(x,{\bf v})=0\right\}
\]
is   an embedded manifold of dimension $d-1$. As a consequence, for every compact subset $K\subset\R^d$, it holds 
\[
\begin{split}
\mathcal{H}^{d-1}\left(\Sigma_{f}\bigcap K\right)&~\leq~\mathcal{H}^{d-1}\left(\left\{\pi(x,{\bf v})\in \R^d:D^{T} f(x){ \bf v}=0,~~(x,{\bf v})\in K\times S^{m-1}\right\}\right)\\
&~\leq~\mathcal{H}^{d-1}\left(\left\{\pi(x,{\bf v})\in \R^d:(x,{\bf v})\in \mathcal{Z}_{\Psi^f}\bigcap  K\times S^{m-1}\right\}\right)~<~\infty,
\end{split}
\]
with $\pi:\R^d\times S^{m-1}$ being the projection on the first component, so that $\pi(x,{\bf v})=x$.
\medskip

{\bf 2.} To complete the proof, we are going to  show that $\mathcal{M}_{\nu}$ is dense in $\C^2(\R^d,\R^m)$  in the next two steps. For any fix $g\in \C^2(\R^d,\R^m)$ and $\ve>0$, we first approximate $g$ by a smooth function $f$ with $\|f-g\|_{\C^2}<\ve$. Then we  need to construct a perturbed smooth function $f^{\theta}$ which is sufficiently close to $f$ in the $\mathcal{C}^2$ norm,  which lies in $\mathcal{M}_{\nu}$. Toward this goal, we define $\varphi:\R^{m^2}\times\R^d\to\R^m$ by
\bel{Theta}
\varphi(\theta,x)~=~\sum_{j=1}^{m}\left(\theta_{jj}+\sum_{j\neq i=1}^m(\theta_{ij}+\theta_{ji})\right) x_j {\bf e}_j,\qquad (\theta,x)\in \R^{m^2}\times \R^d,
\eeq
 with ${\bf e}_1,...,{\bf e}_m$ being the standard basis vectors of $\R^m$ such that 
\bel{con-T}
{\partial \over \partial x_{i}}\varphi_j(\theta,x)~=~\begin{cases}
0&~~~\mathrm{if}\qquad i\in \{m+1,\cdots, d\},\\[2mm]
\theta_{ii}&~~~\mathrm{if}\qquad j=i\in \{1,\cdots, m\},\\[2mm]
\theta_{ij}+\theta_{ji}&~~~\mathrm{if}\qquad j\neq i\in \{1,\cdots, m\}.
\end{cases}
\eeq
For every ${\bf v}\in S^{m-1}$, we compute 
\[
\Psi^{\varphi(\theta,\cdot)}(x,{\bf v})~=~\Pi_{m}\big(D^{T}_{x}\varphi(\theta,x){ \bf v}\big)~=~\sum_{j=1}^{m}\left(\theta_{jj}{\bf v}_j+\sum_{j\neq i=1}^m(\theta_{ij}+\theta_{ji}){\bf v}_i\right) {\bf e}_j,
\]
and  
\[
{\partial\over \partial \theta_{ij}}\Psi^{\varphi(\theta,\cdot)}(x,{\bf v})~=~\begin{cases}
{\bf v}_j{\bf e}_j&~~\mathrm{if}\qquad i=j,\\[2mm]
{\bf v}_i{\bf e}_j+{\bf v}_j{\bf e}_i&~~\mathrm{if}\qquad i\neq j.
\end{cases}
\]
In particular, the matrix $D_{\theta}\Psi^{\varphi(\theta,\cdot)}(x,{\bf v})$ contains the $m\times m$ matrix
\[
P_j~=~\begin{bmatrix}{\bf r}_1\\ \vdots\\ {\bf r}_j\\   \vdots \\{\bf r}_n\end{bmatrix}
~=~
\begin{bmatrix} \bfv_j &0&\cdots &0&\cdots&\cdots&0\cr
0& \bfv_j&\cdots &0&\cdots&\cdots&0\cr
\vdots&\vdots & \ddots&\vdots &\cdots&\cdots&0\cr
{\bf v}_1& \cdots & {\bf v}_{j-1} &{\bf v}_j &{\bf v}_{j+1}&\cdots&{\bf v}_m\cr
0&\cdots&0&0&\bfv_j&\cdots&0\cr
\vdots&&\vdots&\vdots&\vdots&\ddots & \vdots\cr
0&\cdots&0 &0& 0&\cdots &\bfv_j\end{bmatrix},\qquad j\in \{1,\cdots, m\},
\]
and this yields
\bel{rank-T}
\mathrm{rank}\big(D_{\theta}\Psi^{\varphi(\theta,\cdot)}(x,{\bf v})\big)~=~m,\qquad (\theta,x,{\bf v})\in \R^{m^2}\times\R^d\times S^{m-1}.
\eeq

{\bf 3.}  Covering $\overline{B}_{\nu}$ by a finite number of balls $B(x^{\ell},1)$, $\ell=1,2,\cdots, N$, we consider   the family of combined perturbations
\[
f^{\theta}(x)~=~f(x)+\sum_{\ell=1}^{N}\eta\big(x-x^{\ell}\big)\varphi(\theta^{\ell},x-x^{\ell}),\qquad x\in \R^d,
\]
with $\theta=(\theta^1,\cdots, \theta^{N})\in \R^{N d^2}$, $\theta^{\ell}=(\theta^{\ell}_{\ij})\in \R^{d^2}$ and $\eta:\R^d\to [0,1]$ being a smooth cut-off function  such that 
\[
\eta(x)~=~\begin{cases}
1&~~~\mathrm{if}\qquad |x|\leq 1,\\[2mm]
0&~~~\mathrm{if}\qquad |x|\geq 2.
\end{cases}
\]
For every $\ell\in \{1,\cdots, N\}$ and $(x,{\bf v})\in B({\bf x}^{\ell},1)\times S^{m-1}$, it holds
\[
D_{\theta^{\ell}}\Psi^{f^{\theta}}(x,{\bf v})~=~D_{\theta^{\ell}}\Psi^{f}(x,{\bf v})+D_{\theta^{\ell}}\Psi^{\vp(\theta^{\ell},\cdot )}(x-x^{\ell},{\bf v}).
\]
Thus, from (\ref{rank-T}) and $\overline{B}_{\nu}\subseteq \ds\bigcup_{\ell=1}^{N} B(x^{\ell},1)$, we get
\[
\mathrm{rank}\big(D_{\theta}\Psi^{f^{\theta}}(x,{\bf v})\big)~=~m,\qquad (\theta,x,{\bf v})\in \R^{m^2}\times\overline{B}_{\nu}\times S^{m-1}.
\]
By the transversality theorem \cite{RT}, for a dense set of values $\theta\in \R^{m^2}$,  the smooth map $(x,{\bf v})\mapsto \Psi^{f^{\theta}}(x,{\bf v})$    is transversal to $\{0\}$, restricted to the domain $\overline{B}_\nu\times S^{d-1}$. We conclude that the set $\mathcal{M}_{\nu}$ is dense in $\C^2(\R^d,\R^m)$ and the proof is complete. 
\endproof
 
	\section{Quantitative estimates  of critical sets}\label{Sec:3}
	\setcounter{equation}{0}
	In this section, we shall establish a sharp quantitative estimate  on the set of critical points of a  $\mathcal{C}^1$ function $f:[0,1]^{d}\to\R^m$ with $d\geq m$.  More precisely, let's recall  the definition of the set of critical points of $g\in\mathcal{C}^1$,
	\bel{Sm-i}
	\Sigma_g~\doteq~\left\{x\in(0,1)^d:\textrm{rank}(D g(x))~\le~m-1\right\}.
	\eeq
For every  $\ve>0$, we  first establish  an upper bound on  the minimal $(d-1)$-dimensional Hausdorff measure of the sets $\Sigma_g$ over all functions $g\in\C^1([0,1]^d,\R^m)$ within a $\C^1$ distance of $\ve$ from $f$,
	\bel{N-i-f}
	\mathcal{N}_f(\ve)~\doteq~\inf_{g\in B_{\C^1}(f,\ve)} \mathcal{H}^{d-1}\left(\Sigma_g\right).
	\eeq
Before doing this, we shall give  a simple lemma of the decomposition of a unit cube $\square^d=[0,1]^d$ in $\R^d$ which will be used in the proof of Theorem \ref{N-1}.
	\begin{lemma}\label{cub}For every $d\geq 1$, $\square^d$ can be decomposed into $2^{d-1} d!$ polytopes $\Delta_{k}$ in $\R^d$ such that $\Delta_{k}$ has $(d+1)$ vertices  for $k\in\big\{1,2,\dots, 2^{d-1} d!\big\}$.
	\end{lemma}
	{\bf Proof.} The decomposition of $\square^d$ can be done by  using the induction process. Indeed, if $d=1$ then  $\square^1$ is an interval $[0,1]$. Otherwise, for $d\geq 2$,  $\square^d$ has $2d$ faces $\square^{d-1}_h=\partial\square^{d}_{h}$ for $h\in\{1,\dots, 2d\}$ which are  $\R^{d-1}$-cubes of side length $1$. Assume that we can partition each unit cube $\square^{d-1}_h$  into $2^{d-2} (d-1)!$ polytopes $\Delta_{h,k}^{d-1}$ such that $\Delta_{h,k}^{d-1}$ has $d$ vertices  for $k\in\big\{1,2,\dots, 2^{d-2} (d-1)!\big\}$. Then,  $\square^d$ can be partition into $2^{d-1} d!$ polytopes $\Delta_{h,k}^{d}$  for  $h\in\{1,\dots, 2d\}$, $k\in\big\{1,2,\dots, 2^{d-2} (d-1)!\big\}$ such that
		\[
		\Delta_{h,k}^{d}~=~\left\{\theta c+(1-\theta)\cdot y: \theta\in [0,1],y\in \Delta_{h,k}^{d-1}\right\}
		\]
		with $c$ being the center of $\square^d$. 
	\endproof
\subsection{Upper estimate on $\mathcal{N}_f(\ve)$}
	
To obtain an upper  estimate on $\mathcal{N}_f(\ve)$ for a general non-constant $\mathcal{C}^1$-function $f$,  we recall the classical concept of the minimal modulus of  continuity of  continuous functions and their inverse.
	\begin{definition}\label{modu} Given subsets $U\subseteq \R^d$ and $V\subseteq \R^m$, let $h:U\to V$ be continuous. The  minimal modulus of  continuity of $h$ is given by 
		\bel{omega}
	\omega_h(0)~=~0,\qquad \omega_h(\delta)~=~\sup_{x,y\in U,|x-y|\leq \delta} |h(y)-h(x)|\qquad\forall \delta\in (0,\mathrm{diam}(U)].
		\eeq
			The inverse of the minimal modulus of  continuity of   $h$ is the map $[0,\infty)\ni s\mapsto \Psi_h(s)$ defined in $[0,\infty)$ by 
	\bel{Psi}
		\Psi_{h}(s)~:=~\sup\left\{\delta\geq 0: |h(x)-h(y)|\leq s~~\forall |x-y|\leq \delta, x,y\in U\right\}
		\eeq
	\end{definition}
Calling $\omega_{D f}$   the minimal modulus of  continuity of   $D f$, let $\beta_f^{-1}$ be the inverse of the  the strictly increasing map $[0,\infty)\ni\delta\mapsto \beta_{f}(\delta)$ such that
\bel{beta-f}
\beta_{f}(\delta)~=~\omega_{D f}\big(\sqrt{d}\delta\big)\cdot\left[1+\sqrt{d}\delta+4d^{d+\frac{1}{2}}(d+1)^4\right].
\eeq
	\begin{theorem}
		\label{N-1}
		Let $f: [0,1]^d\to\R^m$ be a non-constant  $\mathcal{C}^1$ function with $m\leq d$. Then it holds 
		\begin{equation}
			\label{N-1Eq}
			\begin{split}
				\mathcal{N}_f(\ve)~&\le~d2^{2m+d-1}d!\cdot{1\over \beta_{f}^{-1}(\ve)},\qquad\forall \ve>0.
			\end{split}
		\end{equation}
	\end{theorem}
	{\bf Proof.} The proof is divided into three steps:
	\medskip
	
	{\bf  1.} For a given  $\ve>0$, let $\delta>0$ be a constant such that 
	\bel{delta}
	\beta_f(\delta)~=~\omega_{D f}\big(\sqrt{d}\delta\big)\cdot\left[1+\sqrt{d}\delta+4d^{d+\frac{1}{2}}(d+1)^4\right]~=~\ve.
	\eeq
	We first divide $[0,1]^d$ into $K_\delta\doteq \left\lceil\frac{1}{\delta}\right\rceil^d$ closed cubes $\square_\iota$ of side length $\ell_\delta\le\delta$. For each $\iota \in \{1,\dots, K_{\delta}\}$, we apply Lemma \ref{cub} to partition $\square_\iota$ into $2^{d-1}d!$ polytopes $\Delta_\iota^k$ in $\R^d$ with vertices $v^{k,1}_\iota,\ldots,v_\iota^{k,d+1}$ so that 
	\[
	\Delta_\iota^k~=~\mathrm{co} \left\{v^{k,1}_\iota,\ldots,v_\iota^{k,d+1}\right\}.
	\]
	For every $x\in \Delta_\iota^k$, we write 
	\[
	\begin{split}
		x&~=~\sum_{j=1}^{d+1}\alpha_{\iota}^{k,j}(x)\cdot v^{k,j}_\iota~=~\sum_{j=1}^{d}\alpha_{\iota}^{k,j}(x)\cdot v^{k,j}_\iota+\left(1-\sum_{j=1}^d\alpha_\iota^k(x)\right)\cdot v_\iota^{k,d+1}\\
		&~=~A_\iota^k\cdot\alpha_\iota^k(x)+v_\iota^{k,d+1}
	\end{split}
	\]
	with the $d\times d$ matrix $A_\iota^k$ and vector $\alpha_\iota^k(x)=\big(\alpha_\iota^{k,1},\ldots,\alpha_\iota^{k,d}\big)^{\dagger} \in \R^d$ given by
	\begin{equation}
		\label{alphadef}
		A_\iota^k~=~\left(v_\iota^{k,1}-v_\iota^{k,d+1}\right|\cdots\left|v_\iota^{k,d}-v_\iota^{k,d+1}\right),\qquad \alpha_\iota^k(x)~=~\big(A_\iota^k\big)^{-1}\big(x-v_\iota^{k,d+1}\big).
	\end{equation}
		For each $\iota \in \{1,\dots, K_{\delta}\}$ and $k\in \{1,\dots, 2^{d-1}d!\}$,  the polytope $\Delta_\iota^k$ satisfies  the following geometric properties:
	\begin{itemize}
		\item [(i).] The first $d$ vertices $v_\iota^{k,1},\ldots,v_\iota^{k,d}$ of $\Delta_\iota^k$ are corners of a $d$-dimensional $\ell_\delta$-sided cube $\square_\iota$, and the last vertex $v_\iota^{k,d+1}$ is the center of the cube. Hence, the distance between the first $d$ vertices and the last vertex is
		\begin{equation}
			\label{vdist}
			\left|v_\iota^{k,j}-v_\iota^{k,d+1}\right|~=~\sqrt{\sum_{j=1}^d\left(\frac{\ell_\delta}{2}\right)^2}~=~\frac{\sqrt{d}}{2}\ell_\delta,\qquad  j\in\{1,\ldots,d\},
		\end{equation}
		and this yields 
		\begin{equation}
			\label{diamDelta}
			\text{diam}\big(\Delta_\iota^k\big)~\le~\sqrt{d}\ell_\delta.
		\end{equation}
		\item [(ii).]   The determinant of $A_\iota^k$ is computed by 
		\begin{equation}
			\label{detVol}
			\left|\det\big(A_\iota^k\big)\right|~=~d!\cdot\text{Vol}\big(\Delta_\iota^k\big)~=~{\ell_{\delta}^{d}\over 2^{d-1}}.
		\end{equation}

	\end{itemize}
	{\bf  2.} 
	For every $j\in \{1,\dots, d+1\}$, let $w_\iota^{k,j}:\Delta_\iota^k\to\R^m$ be the  linear approximation of $f$ at the vertex $v_\iota^{k,j}$ such that 
	\begin{equation}
		\label{wdef}
		w_\iota^{k,j}(x)~=~f\big(v_\iota^{k,j}\big)+D f\big(v_\iota^{k,j}\big)\cdot\big(x-v_\iota^{k,j}\big)\quad\forall x\in\Delta_\iota^k.
	\end{equation}
	Set $D_\iota^k(x)=\ds\sum_{\ell=1}^{d+1}\big[\alpha_\iota^{k,\ell}(x)\big]^2$ for every $x\in \Delta_\iota^k$. Recalling the definition of $\alpha^{k}_{\iota}$ in (\ref{alphadef}), we  approximate $f$ in $\Delta_\iota^k$ by
	\begin{equation}
		\label{hiotadef}
		g_\iota^k~=~\sum_{j=1}^{d+1}\beta_\iota^{k,j}\cdot w_\iota^{k,j}\qquad\mathrm{with}\qquad \beta_\iota^{k,j}~\doteq~\ds\frac{\big[\alpha_\iota^{k,j}\big]^2}{D_\iota^k}.
	\end{equation}	
	Then a $\mathcal{C}^1$ approximation $g$ of $f$ in $[0,1]^d$ is defined as follows: For every  $\iota \in \{1,\dots, K_{\delta}\}, k\in \{1,\dots, 2^{d-1}d!\}$, it holds 
	\bel{g-f}
	g(x)~=~g_{\iota}^{k}(x)\qquad\forall x\in \Delta_\iota^k.
	\eeq
	Since   $w_\iota^{k,j}$ are the linear approximations of $f$ at the vertices $v_\iota^{k,j}$,  the value of $g_\iota^{k,j}$ and $D g_\iota^{k,j}$ on a face of $\Delta_\iota^k$ only depends on the values of $f$ and $D f$ at the vertices on that face. In particular, if two polytopes $\Delta_\iota^k$ and $\Delta_{\iota'}^{k'}$ intersect, the approximations $g_\iota^k$ and $g_{\iota'}^{k'}$ and their gradients agree on the intersection. On the other hand, the gradients of $g_\iota^k$ and $\beta^{k,j}_{\iota}$ in $\Delta_\iota^k$ are computed by 
	\begin{equation}
		\label{hderiv1}
		\begin{split}
			D g_\iota^k~&=~\sum_{j=1}^{d+1}w_\iota^{k,j}D\beta_\iota^{k,j}+\beta_\iota^{k,j}\cdot D w_\iota^{k,j},
		\end{split}
	\end{equation}
	\begin{equation}
		\label{betagrad}
		\begin{split}
			\begin{split}
				D\beta_\iota^{k,j}(x)
				&~=~{1\over [D_\iota^k]^2}\cdot \left(\displaystyle2\alpha_\iota^{k,j}\sum_{\ell=1}^{d+1}\alpha_\iota^{k,\ell}\big[\alpha_\iota^{k,\ell}D\alpha_\iota^{k,j}-\alpha_\iota^{k,j}D\alpha_\iota^{k,\ell}\big]\right).
			\end{split}
		\end{split}
	\end{equation}
	From \eqref{hiotadef}, (\ref{diamDelta}), and  (\ref{wdef}), one has   
	\[
	\begin{split}
		|g_\iota^k(x)-f(x)|~&=~\left|\sum_{j=1}^{d+1}\beta_\iota^{k,j}(x)w_\iota^{k,j}(x)-f(x)\right|~=~\left|\sum_{j=1}^{d+1}\beta_\iota^{k,j}(x)w_\iota^{k,j}(x)-\sum_{j=1}^{d+1}\beta_\iota^{k,j}(x)f(x)\right|\\
		&~\le~\max_{1\le j\le d+1}\left|w_\iota^{k,j}(x)-f(x)\right|~\leq~\max_{1\le j\le d+1}\omega_{\nabla f}\Big(\max_{y\in\Delta_\iota^k}\big|v_\iota^{k,j}-y\big|\Big)\cdot\big|x-v_\iota^{k,j}\big|\\
		&~\leq~\omega_{D f}\big(\textrm{diam}\big(\Delta_\iota^k\big)\big)\cdot\textrm{diam}\big(\Delta_\iota^k\big)~\leq~\omega_{D f}\big(\sqrt{d}\delta\big)\cdot\big(\sqrt{d}\delta\big).
	\end{split}
	\]
	Next,  using (\ref{alphadef}) and (\ref{detVol}), we get 
	\[
	\begin{split}
		\max\limits_{j\in\{1,...,d+1\}}\big|D \alpha_\iota^{k,j}(x)\big|&~\le~\sum_{j=1}^d\big|D\alpha_\iota^{k,j}(x)\big|_F~\leq~\sqrt{d}\cdot \big|D\alpha_\iota^{k}(x)\big|_F~=~\sqrt{d}\cdot\big|\big(A_\iota^k\big)^{-1}\big|_F\\
		&~\leq~d\cdot {|A_{\iota}^{k}|^{d-1}_F\over \big|\det\left(A_\iota^k\right)\big|}~\leq~ d\cdot \left({\frac{d\ell_\delta}{2}}\right)^{d-1}\cdot {2^{d-1}\over \ell_{\delta}^{d}}~=~\frac{d^d}{\ell_\delta},
	\end{split}
	\]
	with $|\cdot|_F$ being the Frobenius norm. Thus, from (\ref{betagrad}), we derive   that 
	\bel{grad-b}
	\begin{split}
		\sum_{j=1}^{d+1}\big|D\beta_\iota^k(x)\big|&~\leq~4(d+1)^4\cdot \max\limits_{w\in\{1,...,d+1\}}\big|D\alpha_\iota^{k,j}(x)\big|~\leq~{4d^d(d+1)^4\over \ell_{\delta} }.
	\end{split}
	\eeq
	Notice that $\ds\sum_{j=1}^{d+1}\beta_\iota^{k,j}(x)=1$, $\ds\sum_{j=1}^{d+1}D\beta_\iota^{k,j}(x)=0$, and 
	\[
	\big| w_\iota^{k,j}(x)-f(x)\big|~\leq~\omega_{D f}\left(\sqrt{d}\ell_\delta\right)\cdot\sqrt{d}\ell_\delta,\qquad \big|D w_\iota^{k,j}(x)-D f(x)\big|~\leq~\omega_{D f}\big(\sqrt{d}\delta\big),
	\]
	we estimate
	\[
	\begin{split}
		\left|D g_\iota^k(x)-D f(x)\right|&~=~\left|\sum_{j=1}^{d+1}\beta_\iota^{k,j}(x)D w_\iota^{k,j}(x)+D\beta_\iota^{k,j}(x) w_\iota^{k,j}(x)-D f(x)\right|\\
		&~\leq~\sum_{j=1}^{d+1}\beta_\iota^{k,j}(x)\cdot \big|D w_\iota^{k,j}(x)-D f(x)\big|+\sum_{j=1}^{d+1}\big|D\beta_\iota^{k,j}(x)\big|\cdot \big| w_\iota^{k,j}(x)-f(x)\big| \\
		&~\leq~\omega_{D f}\big(\sqrt{d}\delta\big)\cdot\left(1+4d^{d+\frac{1}{2}}(d+1)^4\right).
	\end{split}
	\]
	Thus, the $\mathcal{C}^1$-distance between $g$ and $f$ is bounded by 
	\bel{C-1-f-g}
	\|g-f\|_{\C^1([0,1]^d,\R^m)}~\le~~\omega_{D f}\big(\sqrt{d}\delta\big)\cdot\left[1+\sqrt{d}\delta+4d^{d+\frac{1}{2}}(d+1)^4\right]~=~\ve.
	\eeq
	
	{\bf  3.} For a fixed $\iota\in  \{1,\dots, K_{\delta}\}$ and $k\in \{1,\dots, 2^{d-1}d!\}$, let's provide an upper bound on the measure of the set 
	\begin{equation}
		\label{sigmaiota}
		\Sigma_{g_\iota^k}~\doteq~\big\{x\in\Delta_\iota^k:\textrm{Rank}\big(D g_\iota^k(x)\big)~\le~m-1\big\}.
	\end{equation}
From (\ref{hiotadef}), (\ref{hderiv1}), and (\ref{betagrad}), $\big[D_{\iota}^{k}(x)\big]^2\cdot Dg_\iota^k(x)$ is a $d\times m$ matrix  of polynomials of degree at most $4$ in $ \Delta_\iota^k$.
In particular,  consider the $d\times d$ matrix
\[
B_\iota^k(x)~=~\left[\big[D_{\iota}^{k}(x)\big]^2\cdot D g_\iota^k(x),{\bf e}_{m+1},...,{\bf e}_d\right],\qquad x\in \Delta_\iota^k,
\]
with ${\bf e}_1,...,{\bf e}_d$ being the standard basis vectors of $\R^d$. Then, its determinant 
\[
p_\iota^k(x)~=~\det\big[B_\iota^k(x)\big],\qquad x\in \Delta_\iota^k,
\]
 is a polynomial of degree at most $2^{2m}$. Notice that $p_\iota^k(x)=0$ for all $x\in \Sigma_{g_\iota^k}$, we have 
		\[
			\Sigma_{g_\iota^k}~\subseteq~\left\{x\in \Delta_\iota^k: p_\iota^k(x)~=~0\right\}.
		\]
Finally, by the main theorem in  \cite{MN0}, we get	
\[
\mathcal{H}^{d-1}\big(\Sigma_{g_\iota^k}\big)~\le~\mathcal{H}^{d-1}\left(x\in \Delta_\iota^k: p_\iota^k(x)~=~0\right)~\leq~d2^{2m}\cdot \delta^{d-1},
\]
and  (\ref{g-f}), (\ref{delta}) yield
		\[
		\begin{split}
		\mathcal{H}^{d-1}\big(\Sigma_{g}\big)~&=~\sum_{\iota=1}^{K_{\delta}}\sum^{2^{d-1}d!}_{k=1}\mathcal{H}^{d-1}\big(\Sigma_{g_\iota^k}\big)~\leq~d{2^{2m+d-1}d!}\cdot {1\over \delta}~=~{d{2^{2m+d-1}d!}\over \beta_{f}^{-1}(\ve)}.
		\end{split}
		\]
The proof is complete.
	\endproof
\medskip
	
We conclude this subsection with the following remark.
	
\begin{remark} For every $\ve>0$ sufficiently small such that $\beta_{f}^{-1}(\ve)\leq \ds{1\over \sqrt{d}}$, the estimate (\ref{N-1Eq}) can be replaced by 
\[
\mathcal{N}_f(\ve)~\le~d^{3/2}{2^{3d-1}d!}\cdot{1\over \Psi_{D f}(\ve/\gamma)}
\]
with $\gamma=2+4d^{d+\frac{1}{2}}(d+1)^4$. In addition, if $f\in C^{1,\alpha}([0,1]^d,\R^m)$ for some $\alpha\in (0,1]$ then 
\[
\Psi_{D f}(s)~\geq~\left({s\over \|D f\|_{\mathcal{C}^{0,\alpha}}}\right)^{1/ \alpha},\qquad s>0,
\]
with $\|D f\|_{\mathcal{C}^{0,\alpha}}=\ds\sup_{x\neq y\in U}{|Df(x)-Df(y)|\over |x-y|^{\alpha}}$ being the  H\"older norm of $D f$. Thus, in this case, it holds 
\[
\mathcal{N}_f(\ve)~\leq~d^{3/2}{2^{3d-1}d!}\cdot\Bigg(\frac{\gamma\|D f\|_{\mathcal{C}^{0,\alpha}}}{\ve}\Bigg)^{1/\alpha}.
\]
\end{remark}

%
\subsection{Lower estimate on $\mathcal{N}_f(\ve)$}		
In this subsection, we show that the upper estimate (\ref{N-1Eq}) is sharp in terms of of power functions.  In order to do so, let's recall the definition of  a modulus of continuity $\omega$ and a function admitting $\omega$ as a modulus of continuity.
	\begin{definition}\label{beta} A function $\omega:[0,\infty]\to[0,\infty]$ is called a modulus of continuity if it is increasing, subadditive, and satisfies
		\[\lim_{\delta\to0+}\omega(\delta)~=~\omega(0)~=~0.\]
		We say that a continuous function $f:U\subset\R^d\to \R^m$ admits $\omega$ as a modulus of continuity if
		\begin{equation}
			\label{betadef}
			\sup_{x,y\in U,|x-y|\le s}|f(x)-f(y)|~\le~ \omega(s)\qquad \forall s\ge0.
		\end{equation}
	\end{definition}	
Our result is stated as follows.
	
	\begin{theorem}\label{lwb}
		Let $d\ge m\ge 1$ and $\omega(\cdot)$ be a modulus of continuity. Then there exists a function $f:[0,1]^d\to\R^m$ such that $D f$ admits $\omega(\cdot)$ as a modulus of continuity and for every $0<\ve<\min\big\{\frac{\beta(2^{-5})}{\gamma},\frac{1}{d}\big\}$, it holds
		\begin{equation}
			\label{lowerBd}
			\mathcal{N}_f(\ve)~\ge~{1\over 16 \Psi_{\omega}(\gamma\ve)}\cdot2^{-2\sqrt{\log_2(\Psi_{\omega}(\gamma\ve))}}
					\end{equation}
with $\gamma=2m^3[\beta(2^{-5})]^{m-1}$ and $\Psi_{\omega}$ being  the inverse of the minimal modulus   of $\omega$.	
	\end{theorem}
{\bf Proof.} {\bf 1.} For every $n\geq 1$, let $c_n:[0,1]\to \R$ be a sample function with $\mathrm{supp}(c_n)\subseteq [0,4\ell_n]$ and $\ell_n=\ds2^{-n^2-n-2}$ such that for every $s\in [0,2\ell_n]$, one has
\[
c_n(s)~=~-c_n(4\ell_n-s)~=~{1\over 2}\cdot\left[\omega(s)\cdot\chi_{[0,\ell_n)}+\omega(2\ell_n-s)\cdot \chi_{[\ell_n,2\ell_n]}(s)\right].
\]
Dividing $[0,1]$ into a countably infinite number of subintervals $[s_n,s_{n+1}]$ with 
\[
s_1~=~0,\qquad s_n=\ds\sum_{\ell=1}^{n}2^{-\ell},\qquad  n\geq 2,
\]
let $\beta:[0,1]\to\R$ be a continuous function such that 
\bel{beta}
\beta(s)~=~\sum_{n=1}^{\infty}\left(\sum_{k=0}^{2^{n^2}-1}c_n(s-s_n-4k\ell_n)\right)\qquad\forall s\in [0,1].
\eeq
Since $c_n(\cdot)$ admits $\omega(\cdot)$ as a modulus of continuity, $\beta(\cdot)$ also  admits $\omega(\cdot)$ as a modulus of continuity. Moreover,  for every $n\geq 1$ and $k\in \{0,1,\cdots, 2^{n^2}-1\}$, the map $s\mapsto \beta$ is decreasing in $[a_{n,k},b_{n,k}]$ and 
 \bel{prb}
 \beta([a_{n,k},b_{n,k}])~=~\big[-\omega(\ell_n)/2,\omega(\ell_n)/2\big]
 \eeq
with
 \[
 a_{n,k}~=~s_n+(4k+1)\ell_n,\qquad  b_{n,k}~=~s_n+(4k+3)\ell_n.
 \]

We claim that the $\mathcal{C}^1$ function $f:[0,1]^d\to\R^m$ defined by 
\[
f(x)~=~\left(\int_0^{x_1}\beta(s)ds,x_2,...,x_m\right),\qquad  x\in [0,1].
\]
is a desired one. 
\medskip

{\bf  2.} The Jacobian matrix $D f(x)\in \R^{m\times d}$ of $f$ which is  is computed by
\[
D f(x)=\begin{bmatrix}\beta(x_1)&0_{1\times (m-1)}&0_{1\times (d-m)}\\0_{(m-1)\times 1}&\mathbb{I}_{m-1}&0_{m-1\times (d-m)}\end{bmatrix},
\]
 admits $\omega(\cdot)$ as a modulus of continuity with  $0_{\ell\times k}\in \R^{d\times k }$ being a zero matrix. Moreover, denote by $\mathbb{I}_{d\times m}=\big[\mathbb{I}_{m\times m}~~ 0_{m\times (d-m)}\big]^{\dagger}$, we have 
 \bel{det1}
 \det(D f(x) \mathbb{I}_{d\times m})~=~ \det\left(\begin{bmatrix}\beta(x_1)&0_{1\times (m-1)}\\0_{(m-1)\times 1}&\mathbb{I}_{m-1}\end{bmatrix}\right)~=~\beta(x_1).
 \eeq
Next we are going show that $f$ satisfies (\ref{lowerBd}) by considering the linear operator $L:\R^{m\times d}\to\R^{m\times m}$ such that 
 \[
 L(A)~=~A \mathbb{I}_{d\times m},\qquad A\in\R^{m\times d}.
 \]
For every $g\in\mathcal{C}^1\big(\R^d,\R^m\big)$ with $\|g-f\|_{\mathcal{C}^1}\leq \ve\leq \frac{\beta(2^{-5})}{2}$, by using Jacobi's formula and the trace properties of a matrix, we estimate 
 \bel{Z-h}
 \begin{split}
\big|\det\big[L(D &g(x))\big]-\det\big[L(D f(x))\big]\big|\\
&=~\left|\int_{0}^1\mathrm{trace}\big(\mathrm{adj}\big(tL(D g(x))+(1-t)L(D f(x))\big) \big[L(D g(x))-L(D f(x))\big]\big)dt\right|\\
&\leq~m\ve\cdot \int_{0}^{1}\left|\mathrm{trace}\big(\mathrm{adj}\big(tL(D f(x))+(1-t)L(D f(g))\big)\right|dt\\
&\leq~m^3\Big[\frac{\beta(2^{-5})}{2}+\ve\Big]^{m-1}\ve~\leq~\big[\beta(2^{-5})\big]^{m-1}m^3\ve
\end{split}
 \eeq
with  $\mathrm{tr}(X)$ and  $\mathrm{adj}(X)$ being the trace and the  adjugate  of the matrix $X$. Moreover, since $\mathrm{rank}(\mathbb{I}_{d\times m})=m$, we have
\[
\mathrm{rank}\big[L(D g(x))\big]~=~\mathrm{rank}[D g(x)],\qquad x\in [0,1],
\]
and this   yields  
\[
\Sigma_g~=~\left\{x\in [0,1]^d:\det\big[L(D g(x))\big]=0\right\}.
\]

{\bf 3.} Finally, to establish a lower bound on $\mathcal{H}^{d-1}(\Sigma_g)$, let $n_0$ be an integer such that 
\bel{n_0}
{\omega(\ell_{n_0+1})\over 2}~\leq~\big[\beta(2^{-5})\big]^{m-1}m^3\ve~\leq~{\omega(\ell_{n_0})\over 2}.
\eeq
Fix $n\in \{1,\cdots, n_0\}$, $k\in \left\{1,2,\cdots, 2^{n^2}-1\right\}$, and ${\bf x}\in [0,1]^{d-1}$. From (\ref{det1})-(\ref{n_0}), one has
\[
\big|\det\big[L(D g({\bf x},s))\big]-\beta(s)\big|~\leq~{\omega(\ell_n)\over 2},\qquad s\in [a_{n,k},b_{n,k}],
\]
and (\ref{prb}) implies that the set 
$$\mathcal{Z}_{n,k}({\bf x})~\doteq~\{s\in [0,1]:\det\big[L(D g({\bf x},s))\big]=0\}$$
is non-empty. In particular, we derive 
\[
\begin{split}
\mathcal{H}^{d-1}\left(\Sigma_g\right)&~\geq~\sum_{n=1}^{n_0}\sum_{k\in \{0,1,\cdots, 2^{n^2}-1\}}\mathcal{H}^{d-1}\left(\bigcup_{{\bf x}\in [0,1]^{d-1}}\mathcal{Z}_{n,k}({\bf x})\times\{{\bf x}\} \right)~\geq~\sum_{n=1}^{n_0}2^{n^2},
\end{split}
\]
and  the first inequality of (\ref{n_0})  yields 
\[
\mathcal{H}^{d-1}\left(\mathcal{Z}_g\right)~\geq~2^{n_0^2}~\geq~{1\over 16 \Psi(\gamma\ve)}\cdot 2^{-2\sqrt{\log_2(\Psi(\gamma\ve))}}.
\]
Thus, (\ref{lowerBd}) holds.
\endproof	
\quad\\
{\small {\bf Acknowledgments.}  This research  was partially supported by National Science Foundation grant DMS-2154201. }

%
%
%

\end{document}